\documentclass{elsarticle}

\usepackage[utf8]{inputenc}
\usepackage[english]{babel}

\usepackage{tikz}
\usetikzlibrary{shapes,arrows}
\usetikzlibrary{calc}

\usepackage{amsfonts,amssymb,amsmath} 
\usepackage{enumerate}

\DeclareMathAlphabet{\mathpzc}{OT1}{pzc}{m}{it}
\newcommand{\cl}{\mathtt{cl}}
\DeclareMathOperator{\interior}{\mathtt{int}}
\newcommand{\Liploc}{\cL\mathpzc{ip}_\mathrm{loc}}
\newcommand{\R}{\mathbb{R}}
\newcommand{\cL}{\mathcal{L}}
\newcommand{\cK}{\mathcal{K}}
\newcommand{\Ras}{\mathbb{R}_{\geq0}}

\newtheorem{defn}{Definition}[section]
\newtheorem{assum}{Assumption}[section]
\newtheorem{thm}{Theorem}[section]
\newtheorem{alg}{Algorithm}[section]

\newtheorem{rem}{Remark}[section]

\journal{Automatica}

\makeatletter
\def\ps@pprintTitle{%
 \let\@oddhead\@empty
 \let\@evenhead\@empty
 \def\@oddfoot{}%
 \let\@evenfoot\@oddfoot}
\makeatother

\begin{document}

\begin{frontmatter}
\title{Almost ISS property for feedback connected systems\tnoteref{mytitlenote}} 
\tnotetext[mytitlenote]{Partially supported by the German Federal Ministry
 of Education and Research (BMBF), research project ''LadeRamProdukt''. Corresponding author P. Feketa.}
\author[Erfurt]{Petro Feketa}\ead{petro.feketa@fh-erfurt.de}
\author[Sydney]{Humberto Stein Shiromoto}\ead{humberto.shiromoto@ieee.org}
\author[Wuerzburg]{Sergey Dashkovskiy}\ead{sergey.dashkovskiy@uni-wuerzburg.de}

\address[Erfurt]{Department of Civil Engineering, University of Applied Sciences Erfurt, Altonaer Str. 25, 99085 Erfurt, Germany}             
\address[Sydney]{The Australian Centre for Field Robotics, The Rose Street Building, J04, The University of Sydney, 2006, NSW, Australia}  
\address[Wuerzburg]{Institute of Mathematics, University of W{\"u}rzburg, Emil-Fischer-Str. 40, 97074 W{\"u}rzburg, Germany}             


\begin{abstract}                          
Small-gain conditions used in analysis of feedback interconnections
are contraction conditions which imply certain stability properties.
Such conditions are applied to a finite or infinite interval. In this paper
we consider the case, when a small-gain condition is applied to several
disjunct intervals and use the density propagation condition in the gaps
between these intervals to derive global stability properties for an interconnection.
This extends and improves recent results from \cite{SteinShiromoto2015}.
\end{abstract}
\end{frontmatter}

\section{Introduction}\label{sec:introduction}
\begin{figure}[width=\textwidth]
\begin{tikzpicture}[
        xscale=2.16,yscale=1.5,
        IS/.style={blue, thick},
        LM/.style={red, thick},
        axis/.style={very thick, ->, >=stealth', line join=miter},
        important line/.style={thick}, dashed line/.style={dashed, thin},
        every node/.style={color=black},
        dot/.style={circle,fill=black,minimum size=4pt,inner sep=0pt,
            outer sep=-1pt},
    ]

    \draw[axis,<->] (2.5,0) node(xline)[right] {$r$} -|
                    (0,2.3) node(yline)[above] {$~$};
    \draw[LM] (0,0) .. controls (1,0.2) and (1.2,1.5) .. (2.3,1.8) node[above] {$\gamma^{-1}_{12}$};
    \draw[IS] (0,0) .. controls (1,0.1) and (1.2,0.7) .. (2.3,1.1) node[below] {$\gamma_{21}$};
    \begin{scope}[xshift=3cm]
        \draw[axis,<->] (0,2.3) node(eyline)[above] {$~$} |-
                        (2.5,0) node(exline)[right] {$r$};

    \draw[IS] plot [smooth] coordinates {(0,0) (0.3,0.4) (0.75,0.55) (1.0,1.3) (1.5,1.5) (2.2,2.1)} node[left] {$\gamma_{21}$};
    \draw[LM] plot [smooth] coordinates {(0,0) (0.3,0.3) (0.7,0.7) (1.3,1.5) (2.2,1.8) } node[below] {$\gamma^{-1}_{12}$};
    \draw[dashed line] (0.45,0.42) node[] {} -- (0.45,0) node[dot, label=below:$\underline M_1$] {} ;
    \draw[dashed line] (0.88,0.93) node[] {} -- (0.88,0) node[dot, label=below:$\overline M_1$] {} ;
    \draw[dashed line] (1.21,1.38) node[] {} -- (1.21,0) node[dot, label=below:$\underline M_2$] {} ;
    \draw[dashed line] (1.73,1.65) node[] {} -- (1.73,0) node[dot, label=below:$\overline M_2$] {} ;
    \end{scope}
\end{tikzpicture}
\caption{Graphs of $\gamma_{21}$ and $\gamma_{12}^{-1}$.}
\label{Fig1}
\end{figure}
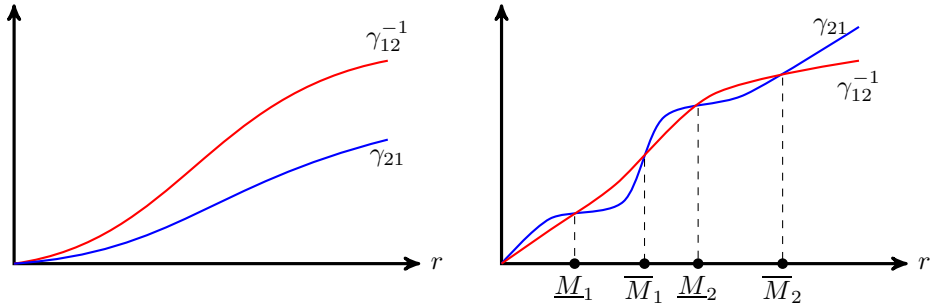

Nonlinear systems are known for essentially more complex behavior compared to linear ones
which makes the investigation of such vital properties as stability and robustness very
challenging and requires essentially different mathematical tools than in case of
linear systems. In many applications nonlinear systems appear in form of interconnections.
In this case the framework of input-to-state stability (ISS) provides useful tools for the stability analysis which includes small-gain theorems \cite{JTP94}
and constructions of the ISS-Lyapunov functions \cite{JIANGETAL:1996}. The small-gain condition requires
that the composition of both interconnection gains satisfies $\gamma_{12}\circ\gamma_{21}(r)<r,\quad r\in(0,\infty)$.
The construction in \cite{JIANGETAL:1996} uses the equivalent statement of the latter condition, namely $\gamma_{21}(r)<\gamma_{12}^{-1}(r),\quad r\in(0,\infty)$, which in particular implies  that the origin is the only common point of the graphs of $\gamma_{21}$ and $\gamma_{12}^{-1}$. See the left part of Fig.~\ref{Fig1}.
In this work we ask: how stability can be studied if the small-gain condition is satisfied not on the whole semi-axis $(0,\infty)$ but on its subset,
or in other words, what can we do in case the graphs of $\gamma_{21}$ and $\gamma_{12}^{-1}$ have several common points (see the right part of Fig.~\ref{Fig1})?
\\
If the small-gain condition is satisfied on an interval of the form $(0,r_0)$ only, then the local ISS property can be established \cite{Dashkovskiy:2010}.
If the small-gain condition fails to hold on some interval $(\underline{M},\overline{M})\subset(0,\infty)$ the work \cite{SteinShiromoto2015}
uses an additional condition to fill this gap and to establish the almost global asymptotic stability of an interconnection.
This condition, called density propagation inequality, is applied to derivatives of functions governing the dynamics
of a system \cite{Rantzer:2001} and is a kind of dual to Lyapunov methods.
\\
In this work we are interested in a more general question: let the small-gain condition be satisfied on some pairwise disjunct intervals
$(\underline{M}_1,\overline{M}_1)\cup\dots\cup(\underline{M}_k,\overline{M}_k)\subset(0,\infty)$, how can we "fill the gaps" between these intervals
by means of suitable density propagation inequalities to assure almost global stability properties for the whole interconnection?
If the points $\underline{M}_1,\overline{M}_1,\dots,\underline{M}_k,\overline{M}_k$ are chosen properly the question of stability can be resolved effectively.
It turns out that these points should be taken as points at which graphs of $\gamma_{21}$ and $\gamma_{12}^{-1}$ coincide. In this work we also provide an algorithm
to find these points.
Our work extends \cite{SteinShiromoto2015} to the case of systems with inputs and multiple regions where the small-gain condition does not hold. Moreover, due to the suitable choice
of these regions we avoid additional restriction (12) in \cite{SteinShiromoto2015}.
This condition requires the inequality $\gamma_{12}^{-1}(\underline{M}) < \gamma_{21}(\overline{M})$ to hold along with the small-gain condition on $(\underline{M}, \overline{M})$.
However not any interval $(\underline{M}, \overline{M})$ where small-gain condition holds satisfies this restriction. Using our new simple algorithm, we can make this interval larger so that (12) is satisfied automatically and hence becomes unnecessary. This is another advantage of this note. An example motivating and illustrating our results is available in \cite{Example}.
\section{Preliminaries and notation}\label{sec:notation}
The notation $\overline{\mathbb{N}}$ (resp. $\overline{\mathbb{R}}$) stands for the set $\mathbb{N}\cup\{\infty\}$ (resp. $\mathbb{R}\cup\{\infty\}$). For a given $a,b\in \overline{\mathbb{R}}$ let $\mathbb N_{[a,b]}=\{s\in\overline{\mathbb{N}}:a\leq s \leq b\}$. Let $\mathbf{S}\subset\R^n$, its closure (resp. interior) is denoted as $\cl\{\mathbf{S}\}$ (resp. $\interior\{\mathbf{S}\}$). We recall the following standard definitions: a function $\alpha:[0,\infty)\to[0,\infty)$ is of class $\mathcal K$ when $\alpha$ is continuous, strictly increasing, and $\alpha(0)=0$. If $\alpha$ is also unbounded, then we say it is of class $\mathcal K_\infty$. A continuous function $\beta:[0,\infty)\times[0,\infty)\to[0,\infty)$ is of class $\mathcal{KL}$, when $\beta(\cdot,t)$ is of class $\mathcal K$ for each fixed $t\geq 0$, and $\beta(r,t)$ decreases to $0$ as $t\to\infty$ for each fixed $r\geq 0$. By $\mathcal{C}^s$ we denote the class of $s$-times continuously differentiable functions, by $\Liploc$ the class of locally Lipschitz continuous functions, by $\mathcal{L}_{\mathrm{loc}}^\infty$ the class of locally essentially bounded functions.

Consider the interconnection of two systems $\Sigma_1$ and $\Sigma_2$
\begin{equation}\label{eq:subsystem}
	\Sigma_i:\quad \dot{x}_i=f_i(x_1(t),x_{2}(t),u_i(t)),\quad i=1,2
\end{equation}
$x_i(t)\in\mathbb{R}^{n_i}$ is the state of $\Sigma_i$ and $u_i(t)\in\mathbb{R}^{m_i}$ is its external input, $f_i$ is assumed to be of class $\mathcal{C}^1$ and satisfy $f_i(0,0,0)=0$. This interconnection can be written as
\begin{equation}\label{newthree}
\dot x=f(x(t),u(t))
\end{equation}
with the state $x=(x_1,x_2)\in\mathbb{R}^n,\; n=n_1+n_2$, dynamics $f=(f_1,f_2)$ and input $u=(u_1,u_2)\in\mathbb{R}^m,\; m=m_1+m_2$.
\begin{defn}\cite{Sontag1989}
	The system \eqref{newthree} is called \emph{input-to-state stable} (ISS) if there exist functions $\beta\in\mathcal{K}\mathcal{L}$ and $\tilde{\gamma}_{u}\in\mathcal{K}_\infty$ such that, for each initial condition $x(0)$ and each measurable essentially bounded input $u(\cdot)$, the solution $x(\cdot)$ to \eqref{newthree} satisfies
	\begin{align*}
		|x(t)|\leq& \beta(|x(0)|,t)+\tilde{\gamma}_{u}(|u|_{\infty}) \quad \forall t\geq 0.
	\end{align*}
\end{defn}
ISS is equivalent to the existence of an ISS-Lyapunov function, which we define for each subsystem in \eqref{eq:subsystem}:
\begin{defn}
	A function $V_i\in\Liploc(\R^{n_i},\Ras)$ is called \emph{storage function} if for some $\underline{\alpha}_i,\overline{\alpha}_i\in\cK_\infty$ it holds that  $\underline{\alpha}_i(|x_i|)\leq V_i(x_i)\leq \overline{\alpha}_i(|x_i|)\quad\forall x_i\in\R^{n_i}$.
\end{defn}
\begin{defn}
	A storage function $V_i$ is called \emph{ISS-Lyapunov function for \eqref{eq:subsystem}} if for some $\gamma_{i,j},\gamma_{i},\alpha_i\in\cK_\infty$ the implication
	\begin{subequations}\label{eq:ISS Lyapunov}
	\begin{equation}\label{eq:ISS Lyapunov:if}
		V_i(x_i)\geq\max\bigg\{\gamma_{i,j}(V_{j}(x_{j})),\gamma_{i}(|u_i|)\bigg\}
	\end{equation}	
	\begin{equation}\label{eq:ISS Lyapunov:then}
	\Rightarrow\;\nabla V_i(x_i)\cdot f_i(x_i,x_{j},u_i)\leq-\alpha_i(|x_i|).
	\end{equation}
	\end{subequations}
  holds for almost all $x_i\in\mathbb R^{n_i}$, $x_j\in\mathbb R^{n_j}$, and $u_i\in\mathbb R^{m_i}$.
$\gamma_{i,j}$ (resp. $\gamma_{i}$) is called \emph{interconnecting} (resp. \emph{external}) \emph{ISS-Lyapunov gain}.
\end{defn}

The stability of the resulting interconnected system \eqref{newthree} can be deduced from the small-gain theorem \cite{JIANGETAL:1996}: if the  interconnecting ISS-Lyapunov gains satisfy
\begin{equation}\label{eq:SGC}
\gamma_{12}\circ\gamma_{21}(s)<s\quad \forall s>0,
\end{equation}
then system \eqref{newthree} is input-to-state stable.

In this paper we assume that the graphs of $\gamma^{-1}_{12}$ and $\gamma_{21}$ have several points of intersection. It means that small-gain condition does not hold globally and the previously known results cannot be utilized to verify global stability properties of the interconnection. To guarantee the desired stability properties of the interconnection, the dual to Lyapunov's techniques \cite{Angeli:2004,Rantzer:2001} is employed in specific domains of the state space.
Our approach extends the results of \cite{SteinShiromoto2015} to the case of arbitrary number of intersection points of $\gamma^{-1}_{12}$ and $\gamma_{21}$ and allows for the presence of non identically zero inputs. Moreover our stability conditions are less restrictive then in \cite{SteinShiromoto2015}.
\section{Main Results}
\begin{assum}[Subsystems are ISS]\label{hyp:ISS}
For each $i=1,2$ there exists an ISS-Lyapunov function $V_i$ for $\Sigma_i$ from \eqref{eq:subsystem} with the corresponding gain functions $\gamma_{i,j},\gamma_{i}\in\cK_\infty$, and $\alpha_i\in\cK_\infty$.
\end{assum}
\begin{assum}[Local SGC]\label{hyp:localSGC}
The small-gain condition $\gamma_{12}\circ\gamma_{21}(s)<s $
holds for any $s\in\mathbf{M}_k=(\underline{M}_k, \overline{M}_{k})$, $k\in\mathbb N_{[1,\ell]}$, $\ell\in{\mathbb{N}}$,
where $\underline{M}_k, \overline{M}_{k}$ correspond to the intersection points of the graphs $\gamma_{12}^{-1}$ and $\gamma_{21}$, i.e., $\gamma_{12}(\gamma_{21}(\underline{M}_k))=\underline{M}_k$ and $\gamma_{12}(\gamma_{21}(\overline{M}_k))=\overline{M}_k$ for any $k\in\mathbb N_{[1,\ell]}$. Here $\overline{M}_l$ can be a finite number or infinity.
\end{assum}
\begin{rem}
The number of intersection points of the graphs of $\gamma_{12}^{-1}$ and $\gamma_{21}$ can be either finite or infinite, with or without accumulation points. It is not assumed that all of them are known, because it is not necessary in the following theorems.
An algorithm to calculate some of these points $\underline{M}_k$, $\overline{M}_{k}$, $k\in\mathbb N_{[1,\ell]}$ is presented in Section~\ref{secalg}. Having a finite number of the intervals as in the above assumption we fill the gaps between them (and at infinity if $\overline{M}_{l}<\infty$) using the density propagation inequality which needs to be satisfied in the corresponding domains of the state space (see Assumption~\ref{hyp:a.e. dissipation inequality} and Theorem~\ref{thm:main result}).
\end{rem}
%
For a given $\delta\in\mathbb R_{\geq 0}$ let $\mathbf{L}_i(\delta)=\{x_i\in \mathbb{R}^{n_i}:V_i(x_i)\leq \delta\}$.
\begin{thm}\label{prop:regional ISS}
Let Assumptions \ref{hyp:ISS} and \ref{hyp:localSGC} hold.
Then there exists $\gamma\in\mathcal K_{\infty}$ such that almost all solutions to system \eqref{newthree} starting in the set $\mathbf{B}_k$ converge to a neighborhood of the set $\mathbf{A}_k$ with radius $\gamma(|u|_{\infty})$, where
\begin{equation}\label{eq:Ak}
\begin{split}
\mathbf{A}_k=\{x\in \mathbb{R}^n:
x_1\in\mathbf{L}_1(\max\{\underline{M}_k, \gamma_{12}(\underline{M}_k)\}),\\ x_2\in\mathbf{L}_2(\max\{\gamma_{21}(\underline{M}_k),\gamma_{21}\circ\gamma_{21}(\underline{M}_k)\})\},
\end{split}
\end{equation}
\begin{equation}\label{eq:Bk}
\mathbf{B}_k=\{x\in \mathbb{R}^n:x_1\in\mathbf{L}_1(\overline{M}_k), x_2\in\mathbf{L}_2(\gamma_{21}(\overline{M}_k))\}.
\end{equation}
If $\overline M_\ell\neq \infty$ then the above mentioned convergence holds only for some bounded inputs $u$.
\end{thm}
{\bf Proof.}
Let us define $V(x)=\max\left\{\sigma(V_1(x_1)), V_2(x_2)\right\}$ for any $x\in\mathbb{R}^n$  and with
\begin{equation}\label{sigma}
\sigma=\frac{\gamma^{-1}_{12}+\gamma_{21}}{2}.
\end{equation}
We follow the ideas of the proof of Theorem~3.1 from \cite{JIANGETAL:1996}. Note that from the Assumption \ref{hyp:localSGC} and \eqref{sigma} it follows that
\begin{equation}\label{help}
\gamma_{21}(r)<\sigma(r)<\gamma_{12}^{-1}(r)\quad \forall r\in\mathbf{M}_k, k\in\mathbb N_{[1,\ell]}.
\end{equation}
Define the following sets:
\begin{equation*}
\begin{split}
A&=\{(x_1,x_2):V_2(x_2)<\sigma(V_1(x_1))\},\\
B&=\{(x_1,x_2):V_2(x_2)>\sigma(V_1(x_1))\},\\
\Gamma&=\{(x_1,x_2):V_2(x_2)=\sigma(V_1(x_1))\}.
\end{split}
\end{equation*}
Now fix any point $p=(p_1,p_2)\neq(0,0)$, and an input value $v=(v_1,v_2)$. There are three cases:

\emph{Case 1: $p\in A$.}
From $p\in A$ follows $V_2(p_2)<\sigma(V_1(p_1))$. From \eqref{help} it follows that $V_1(p_1)>\gamma_{12}(V_2(p_2))$ if $\gamma_{21}(\underline M_k)\leq V_2(p_2)\leq \gamma_{21}(\overline M_k)$. This then implies
\begin{equation}\label{trick}
\nabla V_1(p_1)f_1(p_1,p_2,v_1)\leq -\alpha_1(V_1(p_1))
\end{equation}
whenever $V_1(p_1)\geq \gamma_1(|v_1|)$. It means (see \cite{JIANGETAL:1996}) that for $p\in A_k = A\cap\{(p_1,p_2):\gamma_{21}(\underline M_k)\leq V_2(p_2)\leq \gamma_{21}(\overline M_k)\}$ the following implication holds
\begin{equation*}
V(p)\geq \hat\gamma_1(|v_1|)\quad\Rightarrow\quad \nabla V(p)f(p,v)\leq -\hat\alpha(V(p))
\end{equation*}
 with $\hat\alpha(s)=\sigma'(\sigma^{-1}(s))\alpha_1(\sigma^{-1}(s)),\;
 \hat\gamma_1(r)=\sigma(\gamma_1(r))$.

The cases of $p\in B$ and $p\in \Gamma$ can be treated analogously to \cite{JIANGETAL:1996}. Finally we get that for $p\in A_k\cup B_k\cup \Gamma_k$ there exist $\tilde\alpha, \tilde\gamma \in \mathcal{K}_\infty$ such that
\begin{equation*}
V(p)\geq \tilde\gamma(|v|)\quad\Rightarrow\quad\nabla V(p)f(p,v)\leq -\tilde\alpha(V(p))
\end{equation*}
 with $B_k = B\cap\{(p_1,p_2):V_1(p_1)\in\mathbf{M}_k\}$, $\Gamma_k = \Gamma\cap\{(p_1,p_2):V_1(p_1)\in\mathbf{M}_k\}$.

It remains to check the points that do not fall into the sets $A_k, B_k$, and $\Gamma_k$, $k\in\mathbb N_{[1,\ell]}$. For completeness we put $\underline M_0=\overline M_0 \equiv 0$.
Consider the case $p\in A \cap\{(p_1,p_2):\gamma_{21}(\overline{M}_{k-1})<V_2(p_2)<\gamma_{21}(\underline{M}_k)\}$. For the points $p$ such that $V_1(p_1)>\gamma_{12}(\underline{M}_k)$ it follows that $V_1(p_1)>\gamma_{12}(V_2(x_2))$ which implies \eqref{trick} whenever $V_1(p_1)\geq \gamma_1(|v_1|)$. The rest cases are treated similarly. Combining the previously proved results we conclude that
there exist a storage function $V:\mathbb{R}^n\to\mathbb{R}_{\geq0}$ and functions $\alpha, \gamma\in\mathcal{K}_\infty$ such that, for almost every $x\in\mathbf{B}_k$ (possibly excluding points $\{(x_1,x_2)\in\mathbb R^n: V_1(x_1)=\overline M_k, V_2(x_2)=\gamma_{21}(\overline M_k),k\in\mathbb N_{[1,\ell]}\mbox{~or~}V_1(x_1)=\overline M_\ell\mbox{~or~}V_2(x_2)=\gamma_{21}(\overline M_\ell)\}$), the implication
	\begin{subequations}
		\begin{equation}\label{eq:ISS Lyapunov k:if}
			V(x)\geq\gamma(|u|)\quad\Rightarrow
		\end{equation}
				\begin{equation}\label{eq:ISS Lyapunov k:then}
			\nabla V(x)\cdot f(x,u)\leq-\alpha\left(|x|_{\mathbf{A}_k}\right)
		\end{equation}
holds.
	\end{subequations}	
In the case of $\overline M_\ell\neq\infty$ for sufficiently large inputs the corresponding trajectory may arrive outside the set $\mathbf{B_\ell}$ and then tend to infinity. That is why in the case of $\overline M_\ell\neq \infty$, the convergence stated in Theorem \ref{prop:regional ISS} holds only for bounded inputs $|u|_{\infty}\leq \delta$, $\delta>0$.\hfill$\blacksquare$

The ISS of \eqref{newthree} follows trivially when $\ell=1$, $\underline M_1=0$, $\overline M_1=\infty$; then $\mathbf{A}_1=\{0\}$ and $\mathbf{B}_1= \mathbb{R}^2$. However, when $\ell>1$, solutions to \eqref{newthree} starting in the set $\mathbf{A}_k$ may converge to an $\omega$-limit set  \cite[Birkhoff's Theorem]{Isidori:1995} that lies inside the set $\mathbf{A}_k$ and do not converge to a ball centred at the origin whose radius is proportional to the norm of the input. Due to this fact, the next assumption is needed to check the  asymptotic behaviour of solutions inside the sets $\mathbf{A}_k$. Let $\mathbf{B}_0=\emptyset$ and $\mathbf{A}_{\ell+1}=\mathbb R^n$.

\begin{assum}\label{hyp:a.e. dissipation inequality}
	For each $k\in \mathbb{N}_{[1,\ell+1]}$, there exist an open set
	$\mathbf{D}_k\subset\mathbb{R}^n$ satisfying $\{\mathbf{A}_k\setminus \mathbf{B}_{k-1}\}\subsetneq\cl\{\mathbf{D}_k\}$ and 
	\begin{itemize}
		\item A differentiable function $\rho_k:\mathbf{D}_k\to\mathbb{R}_{>0}$;
		
		\item A continuous function $q_k:\mathbf{D}_k\to\mathbb{R}_{\geq0}$ such that, for almost every $x\in\mathbf{D}_k$, $q_k(x)>0$;
		
		\item A function $\gamma_k\in\mathcal{K}$ such that, for every $x\in\mathbf{D}_k$ and for every $u\in\mathbb R^m$, the following implication holds
		\begin{subequations}
			\begin{equation}\label{eq:}
				\max_{i=1,2}V_i(x_i)\geq\gamma_k(|u|)\quad \Rightarrow
			\end{equation}
			\begin{equation}\label{dpi}
				\div(\rho_k f)(x,u):=\sum_{j=1}^n\tfrac{\partial (\rho_k f_j)}{\partial x_j}(x,u)\geq q_k(x)
			\end{equation}
		\end{subequations}
	\end{itemize}
\end{assum}
\begin{defn}\cite{Angeli:2004}
	The {origin} is called \emph{almost ISS for \eqref{newthree}} if it is locally asymptotically stable and for some $\gamma\in\mathcal{K}_\infty$
	\begin{equation*}
		\limsup_{t\to\infty}|x(t,x(0),u)|\leq\gamma(|u|_\infty)
	\end{equation*}
	holds for every input $u\in\mathcal{L}_{\mathrm{loc}}^\infty(\mathbb{R}_{\geq0},\R^m)$ and for almost every initial condition $x(0)\in \mathbb{R}^n$.
\end{defn}
\begin{thm}\label{thm:main result}
	Under Assumptions \ref{hyp:ISS}, \ref{hyp:localSGC}, and \ref{hyp:a.e. dissipation inequality} system \eqref{newthree} is almost input-to-state stable.
\end{thm}
{\bf Proof.}
Under the imposed assumptions the whole state space is divided into the domains where dissipation implication \eqref{eq:ISS Lyapunov k:if}, \eqref{eq:ISS Lyapunov k:then} holds (due to small-gain condition), and the domains where the density propagation inequality \eqref{dpi} holds. These domains overlap due to Assumption~\ref{hyp:a.e. dissipation inequality}. For any initial condition $x(0)$ there exists $k\in\mathbb N_{[1,\ell]}$ such that either $x(0)\in \mathbf{B}_k$ or $x(0)\in \mathbf{A}_k$. Assumptions \ref{hyp:ISS}, \ref{hyp:localSGC} imply that almost all solutions to system \eqref{newthree} starting in the set $\mathbf{B}_k$ converge to a neighbourhood of the set $\mathbf{A}_k$ with radius proportional to the norm of the input.
Next we show that for almost every initial condition in $\mathbf{A}_k$ and its neighbourhood $\mathbf{D}_k\setminus \mathbf{A}_k$, the corresponding solutions to \eqref{newthree} converge to some neighborhood of the set $\mathbf{B}_{k-1}\setminus\mathbf{D}_k$ with radius proportional to the norm of the input. Let $\mathbf{Z}_k\subset \mathbf{D}_k$
be the set of initial conditions for solutions that remain outside the $\gamma_{k}(|u|_{\infty})$-neighborhood of the set $\mathbf{B}_{k-1}\setminus\mathbf{D}_k$. From \cite{Angeli:2004} and the 2nd part of the proof of Theorem 1 from \cite{SteinShiromoto2015}, Assumption \ref{hyp:a.e. dissipation inequality} implies that the set $\mathbf{Z}_k$ has a Legesgue measure zero. Hence for almost all initial values in $\mathbf{D}_k$ the corresponding solutions converge to a neighborhood of the set $\mathbf{B}_{k-1}\setminus\mathbf{D}_k$ with radius proportional to the norm of the input. Repeating the previous reasoning we conclude the almost ISS of the system \eqref{newthree}.\hfill$\blacksquare$

\begin{rem}
An overlap of 'small-gain' and 'density propagation' regions\linebreak (namely $\{\mathbf{A}_k\setminus \mathbf{B}_{k-1}\}\subsetneq\cl\{\mathbf{D}_k\}$) is essential for the convergence of almost all solutions to a ball centered at the origin with radius $\gamma(|u|_{\infty})$, $\gamma\in \mathcal K_{\infty}$. However properties of this gain function $\gamma$ depend on how large the mentioned overlap is. Too tight overlap, i.e., when the size of the intersection of the mentioned regions is small, leads to a large gain function (i.e. rapidly increasing), and vice
versa, a large overlap can lead to a smaller gain function.
\end{rem}

\section{Algorithm to find gains intersection points}\label{secalg}
This section addresses the question of how to find the intersection of gains $\gamma_{12}^{-1}$ and $\gamma_{21}$ efficiently. If the gains are known we can define $\gamma=\gamma_{12}\circ\gamma_{21}$.
For any $s\geq 0$ one of the following $\gamma(s)< s$, $\gamma(s)> s$ or $\gamma(s)=s$ holds.
In the latter case $s$ is an intersection point we are looking for.	
Let for some $s>0$ be $\gamma(s)< s$ then $\gamma^2(s)=\gamma(\gamma(s))<\gamma(s)<s$, i.e., the sequence $\gamma^n(s),\;n\in\mathbb{N}$ is monotone, decreasing and bounded from below, hence it converges to some $\underline{M}\ge0$. Similarly $(\gamma^{-1})^n$ is either unbounded or converges to some $\overline{M}$. Moreover, it holds that $\gamma(\underline{M})=\underline{M}$ and $\gamma(\overline{M})=\overline{M}$.	

 Let for some $s>0$ be $\gamma(s)> s$, then the sequence $\gamma^n(s),\;n\in\mathbb{N}$ is either unbounded or converges to some $M^\ast$ for which $\gamma(M^\ast)=M^\ast$. The intersection points are those where the small-gain condition is violated.

This leads to the following algorithm. Pick some constant $\Delta\in\mathbb R_{\geq 0}$ (the algorithm's precision). If the sequence of intersection points have an accumulation point, constant $\Delta$ allows to avoid looping the algorithm.
	
	\begin{alg} Input $\Delta>0$. Let $s=0$, $i=1$.
		\begin{enumerate}[{Step} 1.]
\setlength\itemsep{-0.2em}
			\item\label{item:New s} Let $s^\ast=s+\Delta$;
			
			\item If $\gamma(s^\ast)=s^\ast$, go to step \ref{item:Start Gamma(s)=s}. If $\gamma(s^\ast)< s^\ast$, go to step \ref{item:Gamma(s)<s:Start}. If $\gamma(s^\ast)> s^\ast$, go to step \ref{item:Gamma(s)>s:Start};
			
			\item\label{item:Start Gamma(s)=s} Let $s=s^\ast$ and go to step \ref{item:New s};
			
			\item\label{item:Gamma(s)<s:Start} Let $\underline{M}_i=\lim\limits_{n\to\infty}\gamma^n(s^\ast)$ and $\overline{M}_i=\lim\limits_{n\to\infty}\left(\gamma^{-1}\right)^n(s^\ast)$. Go to step \ref{item:Gamma(s)<s:bar M <infty};
			
			\item\label{item:Gamma(s)<s:bar M <infty} If $\overline{M}_i<\infty$ go to step \ref{item:Gamma(s)<s:update s and index}. Else, go to step \ref{item:stop};
			
			\item\label{item:Gamma(s)<s:update s and index} Let $s=\overline{M}_i$ and $i=i+1$, and go to step \ref{item:New s};
			
			\item\label{item:Gamma(s)>s:Start} Let $M^\ast=\lim\limits_{n\to\infty}\gamma^n(s^\ast)$ and go to step \ref{item:Gamma(s)>s:bar M <infty};
			
			\item\label{item:Gamma(s)>s:bar M <infty} If $M^\ast<\infty$, go to step \ref{item:Gamma(s)>s:update s}. Else, go to step \ref{item:stop};
			
			\item\label{item:Gamma(s)>s:update s} Let $s=M^\ast$ and go to step \ref{item:New s};
			
			\item\label{item:stop} $\ell:=i$. Stop algorithm. Leave.
		\end{enumerate}
		Output: $\ell$, and the sets $\{\underline{M}_1,\ldots,\underline{M}_\ell\}$ and $\{\overline{M}_1,\ldots,\overline{M}_\ell\}$.
	\end{alg}

\tikzstyle{decision} = [diamond, draw, fill=green!20,
    text width=4.5em, text badly centered, node distance=2cm, inner sep=0pt]
\tikzstyle{block} = [rectangle, draw, fill=blue!20,
    text width=10em, text centered, rounded corners, minimum height=4em]
\tikzstyle{stopblock} = [rectangle, draw, fill=red!20,
    text width=10em, text centered, rounded corners, minimum height=4em]
\tikzstyle{line} = [draw, -latex']
\tikzstyle{cloud} = [draw, ellipse,fill=red!20, node distance=2cm,
    minimum height=2em]

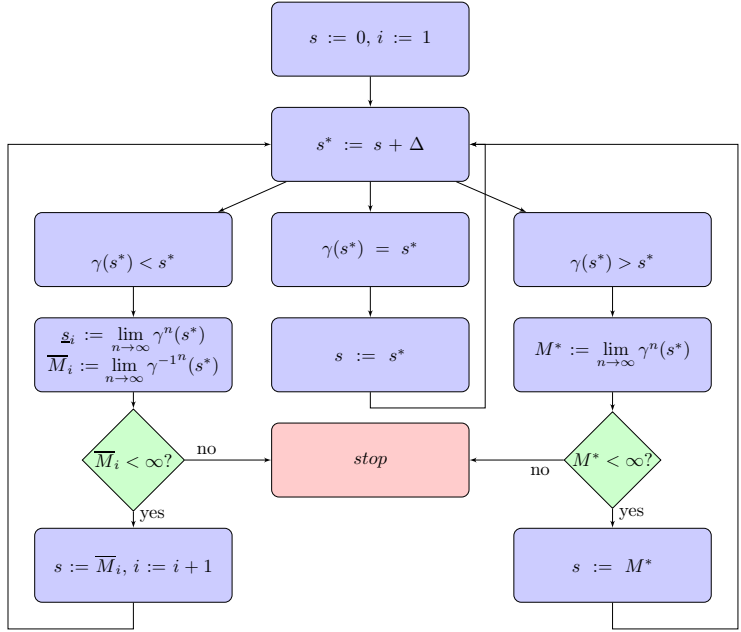
\begin{figure}
\center
 \scalebox{.7}{
\begin{tikzpicture}[node distance = 2cm, auto]
    \node [block] (init) {$s:=0$, $i:=1$};
    \node [block, below of=init] (Delta) {$s^*:=s+\Delta$};
    \node [block, below of=Delta] (equal) {$\gamma(s^*)=s^*$};
    \node [block, left of=equal, node distance=4.5cm] (good) {$$\gamma(s^*)<s^*$$};
    \node [block, right of=equal, node distance=4.6cm] (bad) {$$\gamma(s^*)>s^*$$};

    \node [block, below of=good] (segment) {$\underline s_i:=\lim\limits_{n\to\infty}{\gamma^n(s^*)}$ $\overline M_i:=\lim\limits_{n\to\infty}{{\gamma^{-1}}^n(s^*)}$};
    \node [decision, below of=segment] (decide1) {$\overline M_i<\infty$?};
    \node [block, below of=decide1, node distance=2cm] (reinit1) {$s:=\overline M_i$, $i:=i+1$};

    \node [block, below of=equal] (sss) {$s:=s^*$};

    \node [block, below of=bad] (badsegment) {$M^*:=\lim\limits_{n\to\infty}{\gamma^n(s^*)}$};
    \node [decision, below of=badsegment] (decide2) {$M^*<\infty$?};
    \node [block, below of=decide2, node distance=2cm] (reinit2) {$s:=M^*$};

    \node [stopblock, below of=sss, node distance=2cm] (stop) {$stop$};
    \path [line] (init) -- (Delta);
    \path [line] (Delta) -- (equal);
    \path [line] (Delta) -- (good);
    \path [line] (Delta) -- (bad);
    \path [line] (good) -- (segment);
    \path [line] (segment) -- (decide1);
    \path [line] (equal) -- (sss);
    \path [line] (bad) -- (badsegment);
    \path [line] (badsegment) -- (decide2);

    \path [line] (decide1) -- node [near start] {yes} (reinit1);
    \path [line] (decide1) -- node [near start] {no} (stop);

    \path [line] (decide2) -- node [near start] {yes} (reinit2);
    \path [line] (decide2) -- node [near start] {no} (stop);

    \path [line] (reinit1) |- ($(reinit1.south west) + (-0.5,-0.5)$) |- (Delta);
    \path [line] (reinit2) |- ($(reinit2.south east) + (0.5,-0.5)$) |- (Delta);
    \path [line] (sss) |- ($(sss.south east) + (0.3,-0.3)$) |- ($(Delta.north east) + (0.3,-0.9)$) |- (Delta);

\end{tikzpicture}
}
\caption{A scheme of the numerical-analytic algorithm for small-gain domains allocation}
\end{figure}

The proposed algorithm may have infinitely many iterations (by iteration we mean returning to Step~1) in some cases, as for example in the case of infinitely many points $M_i$ of intersection of gain functions such that $M_i\to\infty$ as $i\to\infty$. Then a reasonable running time should be allowed before termination. A longer running time of the algorithm leads to a smaller domain where the density propagation inequality \eqref{dpi} needs to be satisfied. Hence the algorithm should be interrupted as soon as the density propagation region is small enough to satisfy \eqref{dpi}.

The algorithm ends up in a finite number of iterations for some wide classes of systems: for the systems which gains have either a finite number of intersection points or an infinite number but with finite accumulation points.

\section{Concluding remarks and open problems}
In this paper we have proposed the method of how to verify global stability properties of two feedback connected systems when small-gain condition is not globally satisfied. It was made by imposing additional density propagation inequality in the appropriate state space subsets. The properties of the resulting ISS gain function depend on the overlaps of small-gain and density propagation regions. This dependance is an interesting problem to be investigated in the future.

Another challenging problem is to extend the proposed approach to the case of $n\in\mathbb N$ interconnected systems. This leads to multidimensional small-gain and density propagation regions and related topological and geometrical problems, in particular how to built the corresponding sets $\mathbf{A}_k$ and $\mathbf{B}_k$. A first step to tackle these problems is presented in~\cite{shiromoto2015combination}.

\bibliography{Library}

\end{document}